\documentclass[10pt]{article} 				
\usepackage{amssymb}
\title{Universal Cycles of Complementary Classes}
\author{Michelle Champlin, Anant Godbole \& Beverly Tomlinson\\
East Tennessee State University}
\date{}							% Activate to display a given date or no date

\begin{document}
\maketitle
\abstract{Universal Cycles, or U-cycles, as originally defined by de Bruijn, are an efficient method to exhibit a large class of combinatorial objects in a compressed fashion, and with no repeats.   de Bruijn's theorem states that U-cycles for $n$ letter words on a $k$ letter alphabet exist for all $k$ and $n$.  Much has already been proved about Universal Cycles for a variety of other objects. This work is intended to augment the current research in the area by exhibiting U-cycles for {\it complementary classes}.   Results will be presented that exhibit the existence of U-cycles for class-alternating words such as alternating vowel-consonant (VCVC) words; words with at least one repeated letter (non-injective functions);  words with at least one letter of the alphabet missing (functions that are not onto); words that represent  illegal tournament rankings; and words that do not constitute ``strong" legal computer passwords.  As with previous papers pertaining to U-cycles, connectedness proves to be a nontrivial step.   }

%%%%%%%%%%%%%%%%%%%%%%%%%%%%%%%%%%%%%%%%%%%%%%%%

\section{Introduction}

This article was written as part of an East Tennessee State University Research Experience for Teachers (RET) grant, funded by the National Science Foundation.  The key idea behind the RET program is to expose teachers to open ended problems in mathematics, that they may then use as springboards for open-ended classroom investigations.  The first and third authors intend to do exactly this when they resume their teaching careers in the near future.

Perhaps the most often quoted example of a de Bruijn cycle is 11101000, which, when viewed as a cycle, leads to each of the eight binary three-letter words being represented exactly once.  It is the shortest such representation, as opposed to (for example) the elementary statistics representation of the sample space corresponding to three coin flips, as $$\{111, 110, 101, 100, 011, 010, 001, 000\}.$$  N. G. de Bruijn proved that this is a universal fact that holds true for all length words on any alphabet; thus a U-cycle for all the $k^n$ $n$-letter words on a $k$-letter alphabet exists for all $n,k$.  A key result used in that proof, and, in fact, when verifying existence of U-cycles of several other ``ordered" structures, is the following result from, e.g., \cite{west}.

 \newtheorem{thm} {Theorem}
 \begin{thm}  A digraph is Eulerian if and only if it has at most one non-trivial weakly connected component, and the indegree is equal to the outdegree at every vertex, $i(v)=o(v)$. \end{thm}
One of the results proved in \cite{jackson} was that one could restrict the words to be of a certain kind, and yet show that U-cycles exist; Jackson proved that words in which no letter repeated fell into this category:
\begin{thm} 
Let $k\ge 3$.  Then, a U-cycle of 1-1 functions from $\{1,...n\} \rightarrow \{1,...k\}$ (and of length $(k)_n=k(k-1)\ldots(k-n+1)$) exists if and only if $k>n$; these are merely permutations of $k$ objects taken $n$ at a time, or, $n$-letter words on $[k]$ in which no letter repeats. 
\end{thm}
A result similar to  Theorem 1 was proved in \cite {blg} and presented at the Southeastern Combinatorics Conference in 2008.
\begin{thm} For $k \ge3$, a U-cycle of onto functions from $\{1,...n\} \rightarrow \{1,...k\}$ exists if and only if $n>k$. 
\end{thm}  
Note that when $n=k\ge3$ a U-cycle cannot exist for onto functions {\it or} for one-to-one functions since these both end up being permutations on $[n]$, which break up into disjoint cycles such as $123\rightarrow231\rightarrow312\rightarrow123$.  The objective of this paper is to prove that {\it complementary results} are valid for several of the previously proven results.  Specifically, we will prove existence of U-cycles for words with at least one repeated letter (functions that are {\it not} 1-1),  words with at least one letter of the alphabet missing ({\it non}-onto functions), words that represent  {\it illegal} tournament rankings (the existence of U-cycles of legal rankings was established in \cite{lg}), and words that do {\it not} constitute legal computer passwords (U-cycles of legal passwords were studied in \cite{lg}).  These results will be established in Sections 2 and 3.   Note moreover that {\it a priori} there is no reason for U-cycles of complementary classes to exist; for example it is known \cite {bg} that U-cycles of subsets of $[n]$ of size between $s$ and $t$ exist if $s<t$, but there is no way we can find U-cycles of subsets of size $<s$ and $>t$.  To provide another example, one-to-one functions from $\{A, B, C\}$ to $\{A, B, C\}$ do not admit a U-cycle since these are permutations on [3], but the 21=27-6 non-injective functions {\it can} be arranged in a U-cycle as follows:
$$AAACACCCBBBAABABBCBCC;$$ see Theorem 5 (or Theorem 6) for a general result along these lines.

We end this section  by proving a preliminary (and unrelated, non-complementary) result that exhibits the technique we will continue to use in the later sections.  Our ``class-alternating" Theorem 4  is valid in the context of an alphabet with $k_v$ ``vowels" and $k_c$ ``consonants."  
Thus, $k_v+k_c=k$, where $k$ is our alphabet, and $n$ is the length of the word with alternating vowels and consonants.  We want to prove that ``class alternating words" admit a U-cycle.  Let the vertices of the digraph consist of all class alternating words of length $n-1$, with a directed  edge between two vertices $u,v$ if the last $n-2$ letters of $u$ coincide with the first $n-2$ letters of $v$.  Let the edge label between $u$ and $v$ be given, as is customary in these proofs, by the cocatenation of the vertex labels; it is these edge labels that correspond to the $n$-letter words we wish to form a U-cycle of.    Two cases need to be studied separately: odd $n$ and even $n$.

When $n$ is odd, the words formed along edges both start and end with vowels or both start and end with consonants. A vertex  $v_0=a_1a_2\dots a_{n-1}$ that starts with a vowel and ends with a consonant has outdegree equal to $k_v$ but indegree equal to $k_c$.  This creates the need for a restriction, namely that $k_v=k_c$ so that $i(v_0)=o(v_0)$.  The same restriction is necessary if the vertex starts with a consonant and ends with a vowel.

When $n$ is even the vertices start and end with vowels (or consonants).  In the first case, $i(v_0)=o(v_0)=k_c$.   Alternately, if $a_1$ and $a_{n-1}$ are both consonants, then $i(v_0)=o(v_0)=k_v$.  The (strong) connectedness of the graph is easy to establish.  Let $n$ be odd.  To go from a vertex such as ERACUX to another with the same VC structure such as AMERIC is effortless; to go, e.g., from ERACUX to XANADU takes a little more effort, with an intermediate ``placeholder" step being necessary to achieve the right parity.  We might, for example, start with
ERACUX$\rightarrow$RACUXA, after which the word XANADU can be added on effortlessly in six steps.  The even case is similar.  This establishes that $D=(V,E)$ is Eulerian, and the Euler path gives the required  U-cycle, proving the following result:  
\begin{thm}  There exists a U-cycle of class alternating $n$-letter words on a $k=k_v+k_c$-letter alphabet that consists of $k_v$ vowels and $k_c$ consonants if either $n$ is even, or if $n$ is odd and $k_v=k_c$.
\end{thm} 
NOTE:  Similar results can be formulated and proved if there are three or more categories of words in the alphabet; we skip the details.

\section{Universal Cycles of non-injective 
\newline functions and functions that are not onto}

As in Section 1‍, $k$ and $n$ are defined as the alphabet size and word length through the duration of the paper.  Additionally, the in-degree and out-degree of a vertex are denoted by $i(v)$ and $o(v)$ respectively. 

\begin{thm}
Universal cycles exist for all $n$-letter words on $[k]$ with at least one letter repeated (non-injective functions) provided that $n\ge4$.
\end{thm}

\noindent {\it Proof}: It is true that if $n>k$ {\it all} words are non-injective and thus a U-cycle exists by  de Bruijn's theorem.  However we do not invoke this fact in the proof.  The theorem is false for $n=1$.  If $n=2$, we need $k=1$ for the result to be true.  The proof below relies on having $n\ge4$.  Let $v_1=a_1a_2\dots a_{n-1}$.  There are two cases that we need to consider.  

Case 1: The vertex already has a letter that repeats.

Case 2: Each letter in the vertex is unique and there are no repeats.  

If $v_1$ has repeats, then clearly $i(v_1)=o(v_1)=k$.  If the letters of the word $v_1$ are all unique (this can happen only if $n-1\ge k$), then it is necessary to ``create the repetition." Since there are $n-1$ unique letters in $v_1$, $n-1$ possibilities can create edges with the needed repetition. Therefore, $i(v_1)=o(v_1)=n-1$.

To exhibit (strong) connectivity, define the target vertex $v_2$ to be $b_1b_2\dots b_{n-1}$. If $v_1$ has repeated letters, we combat the possibility that the repeat might be among its first two letters and exhibit the following path between $v_1$ and $v_2$: 

$$\{a_1a_2\dots a_{n-1}\} \rightarrow \{ a_2\dots a_{n-1}b_1\}\rightarrow \{ a_3\dots a_{n-1}b_1b_1\}\rightarrow$$$$ \{ a_4\dots a_{n-1}b_1b_1b_2\}\rightarrow\ldots
\rightarrow \{ b_1b_1b_2\dots b_{n-2}\}\rightarrow \{ b_1b_2\dots b_{n-1}\}.$$

For Case 2, the path can be built as follows:
$$\{a_1a_2\dots a_{n-1}\} \rightarrow \{ a_2\dots a_{n-1}a_{n-1}\}\rightarrow \{ a_3\dots a_{n-1}a_{n-1}b_1\}\rightarrow$$$$ \{ a_4\dots a_{n-1}a_{n-1}b_1b_1\}\rightarrow
\{ a_5\dots a_{n-1}a_{n-1}b_1b_1b_2\}\rightarrow\ldots$$$$\rightarrow \{ b_1b_1b_2\dots b_{n-2}\}\rightarrow \{ b_1b_2\dots b_{n-1}\}$$

In both of the previous cases, a letter is duplicated whenever it is necessary to secure the repetition throughout the mappings.  Therefore, connectivity is proved and a U-cycle exists by Theorem 1. 

\begin{thm}
There exists a universal cycle of functions from $[n]$ to $[k]$ that are not onto, i.e., words that have at least one letter missing, for $n\ge k>2$. \end{thm}

\noindent {\it Proof}:  When $n<k$ all words are non-surjective and thus a U-cycle exists by de Bruijn's theorem, so let us assume that $n\ge k$.   When $k=2$, the result is clearly false since the only possible non-surjective functions are $00\ldots0$ and $11\ldots1$. Consider $v_1=a_1a_2\dots a_{n-1}$, of length $n-1$, from alphabet $[k]$, and with $n\ge k>2$.  The vertex $v_1$ must be restricted to $k-s$  letters, where $s$ represents the number of letters not in the vertex, and $k>s\ge1$.  There are two cases that follow.

Case 1:  $s$ = 1

Case 2: $k>s\ge2$

From the first case, when $s$=1, $v_1$ has $k-1$ unique letters and $i(v_1)=o(v_1)=k-1$, meaning that any letter except the one not in the vertex can be appended to form the edge.

For the second case, when $k>s\ge2$, at least two unique letters from the alphabet are not present. Hence, any letter from the alphabet can be used to form the edge. Thus, $i(v_1)=o(v_1)=k$.

In order to address (strong) connectivity, let the target vertex be $v_2=b_1b_2\dots b_{n-1}$.  $v_1$ maps to $v_2$ smoothly as follows, with $x$ and $y$ denoting any allowable letters that may be appended to $v_1$ and $v_2$ respectively.

$$\{a_1a_2\dots a_{n-1}\} \rightarrow \{ xx\ldots x\}\rightarrow \{ xx\ldots xy\}\rightarrow\{xx\ldots xyy\}\rightarrow\ldots$$$$\rightarrow \{ yy\ldots y\}
\rightarrow \{y\ldots yb_1\}\rightarrow\{y\ldots yb_1b_2\}\rightarrow \{ yb_1\ldots b_{n-2}\}\rightarrow v_2.$$  
Consequently, connectivity is established.  

\section{U-cycles of illegal rankings and words that are not ``strong passwords"}  In this section, we address existence of U-cycles for complementary classes of sets defined in \cite{lg}.  A legal ranking is an ordered word, possibly with repeats, that represent the ranking of contestants in a competition.  For example, 1413 is a legal ranking of the four contestants $a,b,c$ and $d$ -- since it contains no 2, second place being taken up by a tie for first place.  Likewise, 254313 is an illegal ranking.  A legal ranking {\it must} contain a 1.  

\begin{thm}

There exists a universal cycle of illegal rankings of length $n$ (contestants) from an $n$ letter alphabet, where $n\ne 1,3$.    

\end{thm}

\noindent{\it Proof}: First note that the theorem is not true for $n=1$ since the only possible ranking, 1, is legal; the theorem is true for $n=2$ since the only illegal ranking is 22; and the theorem is false for $n=3$ since the three illegal rankings 112, 121, and 211 (out of a total of 27-13=14 illegal rankings) form a connected three cycle.  

Vertices are strings of $n-1$ letters.  These can either be extended to a legal ranking or not.  In the language of \cite{lg} they are either {\it consistent} with a ranking or not.  The in and out degrees $i(v_1)$ and $o(v_1)$ follows a general pattern for all cases. Let $l\ge0$ represent options that can create a legal ranking. Thus, $i(v_1)$ and $o(v_1)$ both are equal to $n-l$, and therefore, $i(v_1)=o(v_1)$.
More specifically, if $v_1$ can be extended to a legal ranking, this is because it is missing a 1, e.g., $v_1=2254$, or another number, e.g. $v_1=1145$.  In the first case $l$ is 1, and in the second case, $l= 2$, since both 11451 and 11453 are legal extensions.  If $v_1$ cannot be extended to a legal ranking, e.g., when $v_1=nn\ldots n$, then $l=0$.

Next, let us examine connectedness.  Let $v_1=a_1a_2 \dots a_{n-1}$.  The target sink vertex that we will show all vertices can lead to (thus establishing weak connectedness) is $v_3=22\ldots2$.  We will reach $v_3$ via an intermediate vertex $v_2$ with no ones.   The mapping of $v_1$ to $v_2$ is as follows:

$$\{a_1a_2\dots a_{n-1}\} \rightarrow \{ a_2\dots a_{n-1}a_n\}\rightarrow\ldots\rightarrow\{a_n\dots a_{2n-2}\}=v_2,$$
where each of the letters $a_n,\ldots,a_{2n-2}$ can be chosen to be different from 1.   Finally, the mapping of $v_2$ to $v_3$ follows with no issues:

$$\{a_n\dots a_{2n-2}\} \rightarrow \{ a_{n+1}\dots a_{2n-2}2\}\rightarrow\ldots\rightarrow \{ a_{2n-2}22\ldots 2\}\rightarrow \{2 \dots 2 2\}.$$

Consequently, since $v_1\rightarrow v_2\rightarrow v_3$, we have weak connectedness, and the result follows from Theorem 1. 
  
\begin{thm}

 Consider a $k$-letter alphabet with $L\ge 3$ categories, and with category $i$ having $k_i$ representatives; $k=\sum_{i=1}^Lk_i$.  Then there exists a U-cucle of $n$-letter words that are non-passwords, i.e., have no letters from at least one category.

\end{thm}
ß†

\noindent{\it Proof}: Strong passwords have all $L$ categories present.  Non-passwords do not meet the criteria.  Note too that the result is false if $L=2$.  Let $v_1$ be a length  $n-1$ word that represents a vertex in the graph.  By definition, $v_1$ has at least one category missing.  Refer to the missing category as $M_1$ when there is exactly one category missing.  In this case, $i(v_1)=o(v_1)=k-\vert M_1\vert$, because $M_1$ is the only category that forms a strong password.  If $v_1$ has two or more categories missing, then clearly $i(v_1)=o(v_1)=k$, because if one of the missing categories is chosen, at least one other category is missing.  

For connectedness, define $v_1=a_1a_2\dots a_{n-1}$ and $v_2=b_1b_2\dots b_{n-1}$. $v_1$ and $v_2$ have at least one category missing, say $M_1$ and $M_2$ respectively, with $M_1=M_2$ being possible.  When considering the mapping we always first choose elements that are not in either of these two categories.  The path from $v_1$ to $v_2$ now follows as indicated, with $B$ being a letter in the complement of the union of categories $M_1$ and $M_2$.

$$\{a_1a_2\dots a_{n-1}\} \rightarrow \{ a_2\dots a_{n-1}B\}\rightarrow \{ a_3\dots a_{n-1}BB\}\rightarrow\dots$$ $$ \rightarrow
\{B\dots B\} \rightarrow \{ B\dots Bb_1\}\rightarrow \dots \rightarrow \{ Bb_1\dots b_{n-2}\}\rightarrow \{ b_1\dots b_{n-1}\}=v_2.$$ 
This proves connectedness.


\begin{thebibliography}{99}
\bibitem{jackson} B. Jackson (1993), ``Universal cycles of $k$-subsets and $k$-permutations,"
{\it Discrete Math.} {\bf 117}, 114--150.
\bibitem{blg} A. Bechel, B. LaBounty-Lay and A. Godbole (2008), ``Universal cycles of
discrete functions," {\it Congressus Numerantium} {\bf 189}, 121--128.
\bibitem{bg} A. Blanca and A. Godbole (2011), ``On Universal cycles for new classes of combinatorial structures," {\it SIAM J. Discrete Math.} {\bf 25}, 1832--1842. 
\bibitem{lg} A. Leitner and A. Godbole (2010), ``Universal cycles of classes of restricted words," {\it Discrete Math.} {\bf 310}, 3303--3309.
\bibitem{west} D. West (1996), Introduction to Graph Theory, Prentice Hall, New Jersey.
\end{thebibliography}
\end{document}